\documentclass{article}

\usepackage{amsmath,
            amssymb,
            latexsym}

\usepackage[pdftex]{graphicx} 
\usepackage{epstopdf}

\newtheorem{theorem}{Theorem}[section]

\newtheorem{corollary}[theorem]{Corollary}
\newtheorem{proposition}[theorem]{Proposition}
\newtheorem{definition}[theorem]{Definition}

\newtheorem{conjecture}[theorem]{Conjecture}

\newenvironment{proof}{\noindent\textsc{Proof: }}
{\hspace{\stretch{1}}$\Box$\medskip}

\begin{document}

\title{A proof of the McKay-Radziszowski\\
subgraph counting conjecture}

\author{
  Alexander Engstr\"om\footnote{The author is Miller Research Fellow 2009-2012 at UC Berkeley, and gratefully acknowledges support from the Adolph C. and Mary Sprague Miller Institute for Basic Research in Science.}
  \\
Department of Mathematics \\
UC Berkeley \\ 
\texttt{alex@math.berkeley.edu}
}

\date\today

\maketitle

\abstract{
We prove a theorem on how to count induced subgraphs in neighborhoods of graphs. Then we use it to
prove a subgraph counting identity conjectured by McKay and Radziszowski in their work on Ramsey theory.
}

\section{Introduction}

This paper is about counting induced subgraphs in neighborhoods. We prove a general theorem that can be used to turn
expressions about enumerations on a normal form: as sums of specific functions. Then we use it to solve a subgraph counting conjecture by McKay and Radziszowski \cite{MR}.

The number of induced subgraphs of $G$ isomorphic to $J$ is called $s(J,G)$. The subgraph of $G$ induced by the
vertices adjacent to $v$ is  $G^{+}_v$, and the one induced by the vertices not adjacent to $v$ is  $G^{-}_v$. This is our main theorem.
\newline
\newline
{\bf Theorem \ref{T2}}\emph{
Let $J_1, J_2, \ldots J_k$ and $J'_1, J'_2, \ldots
J'_l$ be graphs. Then there is a set $\mathcal{J}$ of graphs with at most
\[K=1+\sum_{i=1}^k |V(J_i)| + \sum_{j=1}^l |V(J'_j)| \]
vertices, and constants $m_J$ for each $J\in \mathcal{J}$, such that
\[
 \sum_{v\in V(G)}  
 \prod_{i=1}^k s(J_i,G^-_v)
 \prod_{j=1}^l s(J'_j,G^+_v)
 = 
 \sum_{J\in \mathcal{J}} m_Jj(J,G)
\]
for any graph $G$.
}
\newline
\newline
The function $j(J,G)$ is of inclusion/exclusion type, and is defined by
\[ j(J,G)=\sum_{\phi}  (-1)^{ | (E(\phi(J)) \setminus E(G) | } \]
where $\phi$ ranges over the injections from $V(J)$ into $V(G)$.

In the last section we use Theorem~\ref{T2} to prove this new subgraph counting identity.
\newline
\newline
{\bf Theorem (Conjecture 1 in McKay and Radziszowski \cite{MR})}\emph{
If $G$ is a graph with $n$ vertices, then
\[
\sum_{v\in V(G)}
(p_1(G_v^+)+p_2(G_v^-)+p_3(G_v^+,G_v^-))=0
\]
where
\[\begin{array}{rcl}
p_1(X)&=& n(n-3)s(K_1,X)-(n^2+2n-6)s(K_1,X)^2+3ns(K_1,X)^3 \\
&& -2s(K_1,X)^4+2(n^2+n-8)s(K_2,X)-12s(K_2,X)^2 \\
&& -12(n-1)s(K_1,X)s(K_2,X)
+12s(K_1,X)^2s(K_2,X) \\ && +72s(C_4,X)
    +12(n-2)s(K_3,X)+24s(K_{1,3},X) \\
&& +24s(P_4,X)+24s(T_{3,1},X) 
 +12(n+2)s(P_3,X) \\ && -24s(K_1,X)s(P_3,X)+32s(T_{3,2},X), \\ 
\end{array}\]
\[\begin{array}{rcl}
p_2(Y)&=&4s(K_2,Y)^2-12s(K_{1,3},Y)-8s(C_4,Y)-8s(T_{3,1},Y) 
\\&& -24s(T_{3,2},Y)+2(n-8)s(P_3,Y),
\end{array}\]
and
\[\begin{array}{rcl}
p_3(X,Y)&=&4s(K_1,X)s(P_3,Y)-2(n-2)s(K_1,X)s(K_2,Y)\\ 
&& +4s(K_1,X)^2s(K_2,Y).
\end{array}\]
}

The starting point of this paper was Goodman's identity \cite{G}, and in a sense Theorem~\ref{T2} is a generalization of it.
There are a few infinite families of easily described subgraph counting identities \cite{MR}, and we provide the first example of
one outside them. But there should be an infinite family of "difficult" subgraph counting identities.

\begin{conjecture}
For every $K\geq 4$ there is a subgraph counting identity
\[  \sum_{v\in V(G)}   \sum_{a=1}^N m_a(n) \prod_{i=1}^k s(J_{a,i},G^-_v) \prod_{j=1}^l s(J'_{a,j},G^+_v) = 0 \]
that is true for every graph $G$ with $n$ vertices. In the identity,
%
%
%
\[ \deg m_a(n) + \sum_{i=1}^k |V(J_{a,i})| + \sum_{j=1}^l|V(J'_{a,j})| \leq K \]
for all $a$. There is one graph $J_{a,i}$ or $J'_{a,j}$ with $K$ vertices, and the identity is not in one of the easily described families.
\end{conjecture}

\begin{conjecture}
Modulo the easily described families, and the difficult ones for lower $K$, there is only \emph{one} subgraph counting identity for
every $K\geq 4$.
\end{conjecture}
\section{Counting induced subgraphs}

All graphs and sets are finite. For graph theory not introduced, see Diestel~\cite{D}. 
The complete graph with vertex set $S$ is called $K_S$, and if $S=\{1,2,\ldots,n\}$ it is called $K_n$. If $J$ itself is a graph, then $K_J$ is defined as the
complete graph on the same vertex set as $J$. The graph $T_{n,k}$ is a complete graph $K_n$ extended with an extra vertex that is
adjacent to $k$ of the old ones.

Recall that the function $s(J,G)$ counts the number of induced subgraphs of $G$ isomorphic with $J$.
For any sets $A$ and $B$, the set of injective functions $\phi$ from $A$ to $B$ is
\[ {\tt inj}(A,B). \]
The number of graph automorphisms of $J$ is denoted $\textrm{Aut}(J)$.
Most of our calculations takes place in the polynomial ring
\[ R_S=\mathbb{R}[ x_e \mid e\in E(K_S) ] / \langle x_e^2=1  \mid    e\in E(K_S)  \rangle  . \] 
If $\phi$ is a function defined from the vertex set of a graph, then we use the short $\phi(uv)$ for $\phi(u)\phi(v)$ if $uv$ is an edge. We write $x_{\phi(e)}$ instead
of $x_{\phi(uv)}$ or $x_{\phi(u)\phi(v)}$ if $e$ is the edge $uv$.
Many times the sign of a $\pm 1$ will depend on the truth of a boolean proposition, so we introduce the notation
\[ \pm_P k = \left\{ \begin{array}{rl}
k & \textrm{if $P$ is a true statement,} \\
-k & \textrm{if $P$ is a false statement.}
\end{array}
 \right.  \]

The function $s(J,G)$ enumerates the induced copies of $J$ in $G$. Now we define a polynomial version of it called $\tilde{s}(J,G;x)$. Later on when polynomial versions of functions are introduced, they get the same symbol, but with a tilde on them.

\begin{definition}
Let $J$ be a graph and $S$ a set, then
\[ \tilde{s}(J,G;x) = \frac{
\displaystyle \sum_{\phi\in {\tt inj}(V(J),S)}  \,\, \prod_{e\in K_J}
\left(
1 \pm_{e\in E(J)}x_{\phi(e)}
\right)
}{2^{- |E(K_J)|}  \textrm{\emph{Aut}}(J)  } \]
\end{definition}

The next definition is used for moving between functions and their polynomial versions, as for $s$ and $\tilde{s}$.

\begin{definition}
For any graph $G$ and edge $e\in K_G$, let $x^G_e=\pm_{e\in E(G)}1$.
\end{definition}

\begin{proposition}
For any graphs $J$ and $G$, $s(J,G)=\tilde{s}(J,V(G),x^G)$.
\end{proposition}

\begin{proof}
The product 
\[\prod_{e\in K_J} \left( 1 \pm_{e\in E(J)}x_{\phi(e)}^G \right)\]
is $2^{|E(K_J)|}$ when $\phi(J)=G[\phi(V(J))]$ and zero otherwise. For every induced copy of $J$ in $G$ there are $\textrm{Aut}(J)$ injections from $V(J)$ to $V(G)$
with $\phi(J)=G[\phi(V(J))]$.
\end{proof}

The $\tilde{s}$--polynomials, and products and sums of
them, are broken down into $\tilde{\jmath}$--polynomials.

\begin{definition}
Let $J$ be a graph and $S$ a set. The \emph{$\tilde{\jmath}$--polynomial} in $R_S$ is defined as
\[ \tilde{\jmath}(J,S;x) = \sum_{\phi \in {\tt inj}(V(J),S)} \,\, \prod_{e\in E(J)} x_{\phi(e)} \]
\end{definition}

The corresponding $j$--polynomial is now defined implicitly.
\begin{definition}
Let $J$ and $G$ be graphs, then \[j(J,G)=\tilde{\jmath}(J,V(G);x^G).\]
\end{definition}
There is an easy way to calculate $j$ without $\tilde{\jmath}$.
\begin{proposition}
If $J$ and $G$ are graphs, then
\[ j(J,G)=\sum_{\phi \in {\tt inj}(V(J),V(G)) }  (-1)^{ | (E(\phi(J)) \setminus E(G) | } \]
\end{proposition}

\begin{proof}
\[ \tilde{\jmath}(J,S;x^G) = \sum_{\phi \in {\tt inj}(V(J),S)} \,\, \prod_{e\in E(J)}  \pm_{\phi(e)\in E(G)}1 \]
\end{proof}

\begin{theorem}\label{prop:prodTwoJ}
 Let $J_1$ and $J_2$ be two graphs and $S$ a set. Then
\[ \tilde{\jmath}(J_1,S; x)\tilde{\jmath}(J_2,S; x)=
\sum_{I\subseteq V(J_1)} 
\sum_{\lambda \in {\tt inj}(I, V(J_2))}
\tilde{\jmath}(J_{I,\lambda},S; x)\]
where $J_{I,\lambda}$ is a graph produced like this: Start with the
disjoint union of $J_1$ and $J_2$ and then identify
$v$ in $J_1$ with $\lambda(v)$ in $J_2$ for all $v\in I$.
If there are any double edges after the identifications,
then \emph{both} of the edges should be removed.
\end{theorem}

\begin{proof}
We get the monomials in the sum defining $\tilde{\jmath}(J_k,S;x)$ by relabeling the vertices of $J_k$ with elements of $S$. Given one monomial from $\tilde{\jmath}(J_1,S;x)$
and one from $\tilde{\jmath}(J_2,S;x)$, the $I\subseteq V(J_1)$ on the right hand side of the equality accounts for the vertices of $J_1$ that are relabeled with an
element of $S$ that is also used in the relabeling of $J_2$.

Double edges are removed since $x^2_e=1$ in the ring $R_{S}$ where the $\tilde{\jmath}$--polynomials are defined.
\end{proof}

\begin{corollary}
If $J_1$ is a graph on $n_1$ vertices, and $J_2$ is a graph on $n_2$ vertices,
then $\tilde{\jmath}(J_1,S;x)\tilde{\jmath}(J_2,S;x)$ is a sum of $\tilde{\jmath}(J,S;x)$ polynomials, where no graph $J$ has more than $n_1+n_2$ vertices.
\end{corollary}
Graphs $J$ with isolated points used in $\tilde{\jmath}(J,S;x)$ can be reduced.

\begin{proposition}\label{propIso}
Let $J$ be a graph with an isolated vertex $v$,
and $S$ a set. Then 
\[\tilde{\jmath}(J,S;x)= (|S|-|V(J)|+1)\tilde{\jmath}(J\setminus v,S;x).\]
\end{proposition}
\begin{proof}
The polynomial $\tilde{\jmath}(J,S;x)$ is a sum of monomials indexed by injections from $V(J)$ to $S$. For any injection $\phi$ from $V(J \setminus v)$ to
$S$, we can extend it to an injection from $V(J)$ to $S$ in $(|S|-|V(J)|+1)$ ways without changing the monomial in $\tilde{\jmath}(J,S;x)$.
\end{proof}

The conjecture we prove in this paper is about induced subgraphs inside neighborhoods. And now we have to extend the concept of $\tilde{\jmath}$--polynomials to $\tilde{k}$--polynomials.
For any set $S$, let
\[ R_S^\circ=R_S\cdot
\mathbb{R}[ y_v \mid v \in S ] / \langle y_v^2=1 \mid v\in S \rangle . \] 

\begin{definition} Let $J$ be a graph, $L$ a subset of $V(J)$, and $S$ a set. The \emph{$\tilde{k}$--polynomial} in $R_S^\circ$ is defined as
\[  \tilde{k}(J,L,S;x,y) = \sum_{\phi \in {\tt inj}(V(J),S)}  \,\, \prod_{e\in E(J)} x_{\phi(e)} \prod_{u\in L}y_{\phi(u)}
\]
\end{definition}

The substitution of every $y_u$ in $\tilde{k}(J,L,S;x,y)$ with $x_{uw}$ is denoted
\[ \tilde{k}(J,L,S;x,y)|_{y_u:=x_{uw}}. \]
After the substitution we get a polynomial in $R_{S\cup\{w\}}$.

\begin{proposition}\label{propKJ}
Let $J$ be a graph, $L$ as subset of $V(J)$, and $S$ a set.
Then
\[ \tilde{\jmath}(J',S;x)=\sum_{s\in S} \tilde{k}(J,L,S \setminus s;x,y)|_{y_u:=x_{us}} \]
where $J'$ is the graph $J$ extended by a vertex $w$ with neighborhood $L$.
\end{proposition}

\begin{proof}
\[
\begin{array}{rcl}
\tilde{\jmath}(J',S;x)& =& 
\displaystyle
\sum_{\phi\in {\tt inj}(V(J'),S)}
\,\,
\prod_{e\in E(J')} x_{\phi(e)} \\
& =& 
\displaystyle
\sum_{s\in S}\,\,
\sum_{\phi \in {\tt inj}(V(J'\setminus w),S\setminus s)} \,\,
\prod_{e\in E(J'\setminus w)} x_{\phi(e)} 
\prod_{uw \in E(J')} x_{\phi(u)s} 
\\
& =& 
\displaystyle
\sum_{s\in S}\,\,
\sum_{\phi \in {\tt inj}(V(J),S\setminus s)} \,\,
\prod_{e\in E(J)} x_{\phi(e)}
\bigl.
\prod_{u \in L} y_{\phi(u)} 
\bigr|_{y_u:=x_{us}}
\\
& =& 
\displaystyle
\sum_{s\in S} \tilde{k}(J,L,S
\setminus s;x,y)|_{y_u:=x_{us}}
\end{array}
\]
\end{proof}

Now we define the polynomial equivalent of the $s$--polynomial in the neighborhood case.
\begin{definition}
For any graph $J$ and set $S$, the \emph{$\tilde{r}$--polynomial} in the ring $R_S^\circ$ is defined by
\[  
\tilde{r}(J,S;x,y)= \frac{ \displaystyle
\sum_{V\subseteq V(J)}
\sum_{E\subseteq E(K_J)}
(-1)^{| E\setminus E(J)|}
\tilde{k}(K_J[E],V,S;x,y)
}{2^{-|E(K_J)|-|V(J)|}  \textrm{\emph{Aut}}(J)}
\]
\end{definition}

\begin{proposition}\label{propFacPlus}
If $J$ and $G$ are graphs, then for any vertex $w$ of $G$,
\[  s(J,G_w^+)=\tilde{r}(J,V(G\setminus w);x^G,y)|_{y_u:= x_{uw}^G}.  \]
\end{proposition}

\begin{proof}
By definition $s(J,G^+_w)$ equals $\tilde{s}(J,G[N(w)];x^G)$. 
We reformulate the
product
$2^{|E(K_J)|}  \textrm{Aut}(J)   ^{-1}\tilde{s}(J,G[N(w)];x^G)$ several times:
\[  \sum_{\phi \in {\tt inj}(V(J),N(w))} \prod_{e\in E(K_J)}(1\pm_{e\in E(J)}x_{\phi(e)}) \]
\[  \sum_{\phi \in {\tt inj}(V(J),V(G\setminus w))} 
 \biggl( \prod_{v\in V(J)} \frac{ 1+ x^G_{\phi(v)w}}{2} \biggr)
 \biggl( \prod_{e\in E(K_J)}(1\pm_{e\in E(J)}x_{\phi(e)}) \biggr) \]
 \[   2^{-|V(J)|} \makebox[-1cm]{}
 \sum_{\phi \in {\tt inj}(V(J),V(G\setminus w))}
 \makebox[-0.15cm]{}
  \biggl(
\sum_{V\subseteq V(J)}
\prod_{v\in V} x^G_{\phi(v)w}
 \makebox[-0.1cm]{}
\biggr)
\biggl(
\sum_{E \subseteq E(K_J) }
(-1)^{|E\setminus E(J)|}
\prod_{ e \in E} 
x_{\phi(e)}
 \makebox[-0.1cm]{}
\biggr)\]
\[
2^{-|V(J)|}   \makebox[-0.25cm]{}
\sum_{V\subseteq V(J)}
\sum_{E \subseteq E(K_J) }
(-1)^{|E\setminus E(J)|} \makebox[-0.8cm]{}
\sum_{\phi \in {\tt inj}(V(J),V(G\setminus w))}  \,
\prod_{v\in V} x^G_{\phi(v)w}
\prod_{ e \in E} 
x_{\phi(e)}
\]
\[
2^{-|V(J)|}  \makebox[-0.25cm]{}
\sum_{V\subseteq V(J)}
\sum_{E \subseteq E(K_J) }
(-1)^{|E\setminus E(J)|}
\tilde{k}
(K_J[E],V,V(G\setminus w);x^G,y)|_{y_u:=x^G_{uw}}
\]
and get the desired
 $2^{|E(K_J)|}  \textrm{Aut}(J)^{-1}\tilde{r}(J,V(G\setminus w);x^G,y)|_{y_u:= x_{uw}^G}.$
\end{proof}

We also need a polynomial equivalent of $s$ for when we are outside the neighborhood.

\begin{definition}
For any graph $J$ and set $S$, the \emph{$\tilde{q}$--polynomial} in the\linebreak ring $R_S^\circ$ is defined by
\[  
\tilde{q}(J,S;x,y)= \frac{ \displaystyle
\sum_{V\subseteq V(J)}
\sum_{E\subseteq E(K_J)}
(-1)^{|V|+ |E\setminus E(J)|}
\tilde{k}(K_J[E],V,S;x,y)
}{2^{-|E(K_J)|-|V(J)|}     \textrm{\emph{Aut}}(J)}
\]
\end{definition}

\begin{proposition}\label{propFacMinus}
If $J$ and $G$ are graphs, then for any vertex $w$ of $G$,
\[s(J,G^-_w)=\tilde{q}(J,V(G\setminus w ),
              x^G,y)|_{y_u=x^G_{uw}}.\]
\end{proposition}

\begin{proof}
The proof of Proposition~\ref{propFacPlus} goes through with minor modifications. The calculation is done
for the vertices outside the neighborhood of $w$, so the term
\[ \prod_{v\in V(J)} \frac{1+x^G_{\phi(v)w}}{2} \]
is changed to
\[ \prod_{v\in V(J)} \frac{1-x^G_{\phi(v)w}}{2} \]
and that contributes the extra $(-1)^{|V|}$ in the definition of $\tilde{q}$ compared to $\tilde{r}$.
\end{proof}

\begin{theorem}\label{propMulK}
 Let $J_1$ and $J_2$ be two graphs,
 $L_1\subseteq V(L_1), L_2\subseteq V(L_2)$, and $S$ a set. 
Then
$ \tilde{k}(J_1,L_1,S;x,y)\tilde{k}(J_2,L_2,S;x,y)$ equals
\[
\sum_{I\subseteq V(J_1)} \,\,
\sum_{\lambda \in {\tt inj}(I,V(J_2)) }
\tilde{k}(J_{I,\lambda},L_{I,\lambda},S,x,y)\]
where $J_{I,\lambda}$ is a graph and $L_{I,\lambda}$ a subset of
its vertex set produced like this: Start with the
disjoint union of $J_1$ and $J_2$ where the vertices
of $J_i$ that are in $L_i$ are colored blue and the
other ones red.
Then identify
$v$ in $J_1$ with $\lambda(v)$ in $J_2$ for all $v\in I$.
If $v$ and $\lambda(v)$ had the same color before
their identification the identified vertex is red
and otherwise blue. If there are any double edges 
after the identifications then \emph{both} of the 
edges should be removed. The graph we get is $J_{I,\lambda}$
and the blue vertices form $L_{I,\lambda}$.
\end{theorem}

\begin{proof}
The proof is a straightforward generalization
of the proof of Theorem~\ref{prop:prodTwoJ} and
is left to the reader. 
\end{proof}

\begin{corollary}\label{cor:MulK}
Let $J_1$ and $J_2$ be two graphs,
 $L_1\subseteq V(L_1), L_2\subseteq V(L_2)$, and $S$ a set. 
Then
$ \tilde{k}(J_1,L_1,S;x,y)\tilde{k}(J_2,L_2,S;x,y)$ equals a sum of\linebreak
$ \tilde{k}(J,L,S;x,y)$ polynomials where the graph $J$ has at most
$|V(J_1)|+|V(J_2)|$ vertices.
\end{corollary}

This is our main theorem.

\begin{theorem}\label{T2}
Let $J_1, J_2, \ldots J_k$ and $J'_1, J'_2, \ldots
J'_l$ be graphs. Then there is a set $\mathcal{J}$ of graphs with at most
\[K=1+\sum_{i=1}^k |V(J_i)| + \sum_{j=1}^l |V(J'_j)| \]
vertices, and constants $m_J$ for each $J\in \mathcal{J}$, such that
\[
 \sum_{v\in V(G)}  
 \prod_{i=1}^k s(J_i,G^-_v)
 \prod_{j=1}^l s(J'_j,G^+_v)
 = 
 \sum_{J\in \mathcal{J}} m_Jj(J,G)
\]
for any graph $G$.
\end{theorem}

\begin{proof}
Using Proposition~\ref{propFacMinus} any $s(J_i,G^{-}_v)$ can be expanded into the form
\[ \sum_{V\subseteq V(J_i)} \sum_{E \subseteq E(K_{J_i})}
k_{i,V,E} \tilde{k}(K_{J_i}[E],V,V(G\setminus v);x^G,y)|_{y_u:=x^G_{uv}}
\]
and any $s(J_i',G^{+}_v)$ can be turned into 
\[ \sum_{V\subseteq V(J_i)} \sum_{E \subseteq E(K_{J_i})}
k_{i,V,E}' \tilde{k}(K_{J_i'}[E],V,V(G\setminus v);x^G,y)|_{y_u:=x^G_{uv}}
\]
by Proposition~\ref{propFacPlus}.

Observe that the coefficients $k_{i,V,E}$ and $k'_{j,V,E}$ does \emph{not} depend on $v$ or $G$. By Corollary~\ref{cor:MulK} the product of  $\tilde{k}(J,\ldots)$ and $\tilde{k}(J',\ldots)$ is a weighted sum of $\tilde{k}(J'',\ldots)$ polynomials
with $|V(J'')| \leq |V(J)|+|V(J')|$. 
Thus the products
\[\prod_{i=1}^k  \sum_{V\subseteq V(J_i)} \sum_{E \subseteq E(K_{J_i})}
k_{i,V,E} \tilde{k}(K_{J_i}[E],V,V(G\setminus v);x^G,y)
\]
and
\[\prod_{i=1}^l \sum_{V\subseteq V(J_i)} \sum_{E \subseteq E(K_{J_i})}
k_{i,V,E}' \tilde{k}(K_{J_i'}[E],V,V(G\setminus v);x^G,y)
\]
can be rewritten as
\[  \sum_{J\in \mathcal{J}'}
  \sum_{L\in V(J)}
   m_{J,L}\tilde{k}(J,L,V(G\setminus v);x^G,y) \]
where $m_{J,L}$ are independent of $G$ and $v$; and $\mathcal{J}'$
is a set of graphs with at most $K-1$ vertices.
By Proposition~\ref{propKJ}
\[ \sum_{v\in V(G)}  \sum_{J\in \mathcal{J}'} \sum_{L\in V(J)}
    m_{J,L}\tilde{k}(J,L,V(G\setminus v);x^G,y)|_{y_u=x^G_{uv}}
\]
can be written as a sum  
\[\sum_{J\in \mathcal{J}} m_J\tilde{\jmath}(J,V(G);x^G) = 
\sum_{J\in \mathcal{J}} m_Jj(J,G)  \]
where $\mathcal{J}$
is a finite set of graphs on at most $K$ vertices.
\end{proof}

We do not want graphs in $\mathcal{J}$ with isolated vertices.

\begin{corollary}\label{workie}
Let $J_1, J_2, \ldots J_k$ and $J'_1, J'_2, \ldots
J'_l$ be graphs. 
Then there is a finite
set $\mathcal{J}$ of graphs without isolated vertices, and
polynomials $m_J(n)$ such that
\[ |V(J)| + \deg m_J(n) \leq 
1+\sum_{i=1}^k |V(J_i)| + \sum_{j=1}^l |V(J'_j)| \]
and
\[
 \sum_{v\in V(G)}  
 \prod_{i=1}^k s(J_i,G^-_v)
 \prod_{j=1}^l s(J'_j,G^+_v)
 =
 \sum_{J\in \mathcal{J}} m_J(n)j(J,G)
\]
for any graph $G$.
\end{corollary}

\begin{proof}
For every isolated vertex removed with Proposition~\ref{propIso}
the degree of the corresponing polynomial increases with
at most one.
\end{proof}

\section{A proof of the McKay-Radziszowski\\ subgraph counting conjecture}

In this section we prove the McKay-Radziszowski subgraph counting conjecture stated in the introduction. 

\begin{proof}
In the conjectured equality there are many terms of the type
\[  \sum_{v\in V(G)}  
 \prod_{i=1}^k s(J_i,G^-_v)
 \prod_{j=1}^l s(J'_j,G^+_v).
\]
Our proof strategy is to expand all of them using Corollary~\ref{workie}, and then show that the terms cancel out. 

To calculate that the terms cancels, given the expansions, is elementary but tedious. We have performed it both by hand and by computer, and will not write down the calculations in this paper.

We will however tabulate the different $m_J(n)$ polynomials whose existence is given by Corollary~\ref{workie}. We calculated them by inserting all graphs $G$ with less than 10 vertices, and many large random graphs,  to get very determined linear equations for the coefficients in the $m_J(n)$ polynomials.

The graphs on less than six vertices with no
isolated vertices are 
\[
\begin{array}{ccc}
G_2=K_2   	& G_3=P_3   	& G_4=K_3\\
G_5=K_{1,3}   	& G_6=2K_2  	& G_7=P_4\\ 
G_8=T_{3,1}  	& G_9=C_4 	& G_{10}=K_4-K_2  \\ 
G_{11}=K_4    	& G_{12}=K_{1,4}     & G_{13}=K_2\cup P_3\\ 
G_{15}=K_{1,4}+K_2   & G_{16}=P_5     & G_{20}=K_3\cup K_2 \\ 
G_{23}=T_{4,1}   & G_{24}=K_{2,3}       & G_{25}=K_5-K_3  \\ 
G_{27}=C_5           & G_{28}=C_5+K_2  & G_{29}=K_5-P_4  \\   
G_{30}=K_5-(P_3\cup K_2)          & G_{31}=K_5-P_3  &  G_{32}=K_5-2K_2 \\ 
G_{33}=K_5-K_2  &  G_{34}=K_5 \\
\end{array}\]
and the rest of them are drawn in Figure~\ref{fig}.
\begin{figure}
  \begin{center}
  \includegraphics*{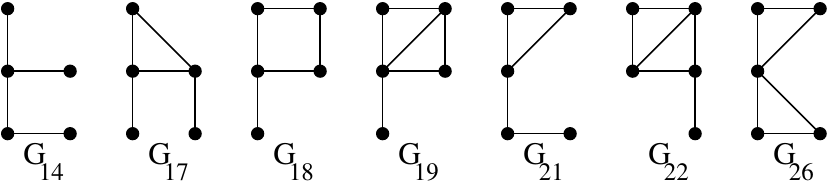}
  \caption{Some graphs on 5 vertices}\label{fig}
  \end{center}
 \end{figure}
To shorten the list we use the notation $n^ {\underline m} = n(n-1)(n-2)\cdots (n-m+1).$
Let $G$ be a graph with $n$ vertices, and let $g_i=j(G_i,G)$ for $i=2,3,\ldots, 34$. To save space, $\sum=\sum_{v\in V(G)}$ in the table. The expansions from Corollary~\ref{workie} are as follows: 
\newline
\newline
$2 \sum s(K_1,G_v^+)=n^{\underline 2}+g_{2}$,
\newline
\newline
$2^{2} \sum s(K_1,G_v^+)^2=nn^{\underline 2}+2(n-1)g_{2}+g_{3}$,
\newline
\newline
$2^{3} \sum s(K_1,G_v^+)^3=n^{\underline 2}({n}^{2}+n-2)+(3{n}^{2}-3n-2)g_{2}+3(n-1)g_{3}+g_{5}$,
\newline
\newline
$2^{4} \sum s(K_1,G_v^+)^4=nn^{\underline 2}({n}^{2}+3n-6)+4(n-1)({n}^{2}+n-4)g_{2}+2(3{n}^{2}-3n-4)g_{3}+4(n-1)g_{5}+g_{12}$,
\newline
\newline
$2^{4} \sum s(K_2,G_v^+)=n^{\underline 3}+3(n-2)g_{2}+3g_{3}+g_{4}$,
\newline
\newline
$2^{8} \sum s(K_2,G_v^+)^2=n^{\underline 3}({n}^{2}+n+4)+2n(n-2)(3n-1)g_{2}+2(5{n}^{2}-3n-12)g_{3}+2({n}^{2}+n-4)g_{4}+4ng_{5}+(5n-4)g_{6}+8(n+2)g_{7}+4(n+4)g_{8}+8g_{9}+8g_{10}+g_{12}+6g_{13}+4g_{14}+2g_{15}+4g_{16}+2g_{20}+4g_{21}+g_{26}$,
\newline
\newline
$2^{5} \sum s(K_1,G_v^+)s(K_2,G_v^+)=n^{\underline 3}(n+1)+4n(n-2)g_{2}+(5n-3)g_{3}+(n+1)g_{4}+g_{5}+g_{6}+2g_{7}+g_{8}$,
\newline
\newline
$2^{6} \sum s(K_1,G_v^+)^2 s(K_2,G_v^+)=n^{\underline 3}({n}^{2}+3n-2)+(n-2)(5{n}^{2}+5n-12)g_{2}+2(4{n}^{2}-3n-3)g_{3}+({n}^{2}+3n-2)g_{4}+2(2n-3)g_{5}+2(n+1)g_{6}+4(n+1)g_{7}+2(n+1)g_{8}+g_{12}+g_{13}+2g_{14}+g_{15}$,
\newline
\newline
$2^{13} \sum s(C_4,G_v^+)=n^{\underline 5}+6(n-2)^{\underline 3}g_{2}+6(n-3)^{\underline 2}g_{3}-2(n-3)^{\underline 2}g_{4}-4(n-4)g_{5}+7(n-4)g_{6}+4(n-4)g_{7}-20(n-4)g_{8}-(n-4)g_{9}-14(n-4)g_{10}-3(n-4)g_{11}-3g_{12}+10g_{13}-4g_{14}-14g_{15}+12g_{16}-12g_{17}+4g_{18}-20g_{19}+2g_{20}+4g_{21}-4g_{22}-4g_{23}+2g_{24}-2g_{25}-g_{26}+4g_{27}+12g_{28}+4g_{29}+10g_{30}+6g_{31}+7g_{32}+6g_{33}+g_{34}$,
\newline
\newline
$3\cdot 2^{7} \sum s(K_3,G_v^+ )=n^{\underline 4}+6(n-2)^{\underline 2}g_{2}+12(n-3)g_{3}+4(n-3)g_{4}+4g_{5}+3g_{6}+12g_{7}+12g_{8}+3g_{9}+6g_{10}+g_{11}$,
\newline
\newline
$3\cdot 2^{11} \sum s(K_{1,3},G_v^+)=n^{\underline 5}+4(n-2)^{\underline 3}g_{2}+6(n-3)^{\underline 2}g_{3}-4(n-3)^{\underline 2}g_{4}+8(n-4)g_{5}-3(n-4)g_{6}-12(n-4)g_{8}+3(n-4)g_{9}-12(n-4)g_{10}-5(n-4)g_{11}+5g_{12}-12g_{13}+12g_{14}+12g_{15}-12g_{16}+24g_{18}+12g_{19}-10g_{20}-24g_{21}-12g_{22}-8g_{23}+10g_{24}+4g_{25}-3g_{26}+12g_{28}+12g_{30}-6g_{31}+3g_{32}-4g_{33}-g_{34}$,
\newline
\newline
$2^{11} \sum s(P_4,G_v^+)=n^{\underline 5}+4(n-2)^{\underline 3}g_{2}+2(n-3)^{\underline 2}g_{3}+(n-4)g_{6}-8(n-4)g_{7}-4(n-4)g_{8}-5(n-4)g_{9}-4(n-4)g_{10}-(n-4)g_{11}+g_{12}-4g_{14}+4g_{15}-4g_{16}-8g_{18}+4g_{19}+2g_{20}+8g_{21}+4g_{22}-2g_{24}+5g_{26}+4g_{28}+8g_{29}-2g_{31}-g_{32}-4g_{33}-g_{34}$,
\newline
\newline
$2^{11} \sum s(T_{3,1},G_v^+)=n^{\underline 5}+6(n-2)^{\underline 3}g_{2}+10(n-3)^{\underline 2}g_{3}+2(n-3)^{\underline 2}g_{4}+4(n-4)g_{5}+3(n-4)g_{6}+4(n-4)g_{7}+4(n-4)g_{8}-(n-4)g_{9}+2(n-4)g_{10}+(n-4)g_{11}+g_{12}-2g_{13}-4g_{14}+2g_{15}-12g_{16}-4g_{17}-12g_{18}+4g_{19}-2g_{20}-12g_{21}-4g_{22}+4g_{23}-2g_{24}+2g_{25}-g_{26}-4g_{27}-12g_{28}+4g_{29}-2g_{30}+10g_{31}+3g_{32}+6g_{33}+g_{34}$,
\newline
\newline
$2^{7} \sum s(P_3,G_v^+)=n^{\underline 4}+4(n-2)^{\underline 2}g_{2}+4(n-3)g_{3}+g_{6}-4g_{8}-g_{9}-4g_{10}-g_{11}$,
\newline
\newline
$2^{8} \sum s(K_1,G_v^+)s(P_3,G_v^+)=n^{\underline 4}(n+2)+(n-2)^{\underline 2}(5n+4)g_{2}+(n-3)(7n-4)g_{3}+3(n-4)g_{5}+2(n-1)g_{6}+2(n-4)g_{7}-3(n+4)g_{8}-(n+2)g_{9}-4(n+2)g_{10}-(n+2)g_{11}+g_{12}+g_{14}+g_{15}-2g_{16}-2g_{17}-g_{18}-g_{19}-g_{20}-3g_{21}-3g_{22}-g_{23}$,
\newline
\newline
$2^{12} \sum s(T_{3,2},G_v^+)=n^{\underline 5}+8(n-2)^{\underline 3}g_{2}+18(n-3)^{\underline 2}g_{3}+4(n-3)^{\underline 2}g_{4}+8(n-4)g_{5}+9(n-4)g_{6}+24(n-4)g_{7}+12(n-4)g_{8}+3(n-4)g_{9}-(n-4)g_{11}+g_{12}+12g_{13}+12g_{14}+12g_{16}-12g_{19}+2g_{20}-12g_{22}-8g_{23}-2g_{24}-4g_{25}-3g_{26}-12g_{28}-24g_{29}-12g_{30}-18g_{31}-9g_{32}-8g_{33}-g_{34}$,
\newline
\newline
$2^{8} \sum s(K_2,G_v^-)^2=n^{\underline 3}({n}^{2}+n+4)-2(n-2)({n}^{2}-3n+8)g_{2}+2({n}^{2}-7n+4)g_{3}+2({n}^{2}+n-4)g_{4}-4ng_{5}-(3n+4)g_{6}+8(n-2)g_{7}-4(n-4)g_{8}+8g_{9}-8g_{10}+g_{12}-2g_{13}-4g_{14}+2g_{15}+4g_{16}+2g_{20}-4g_{21}+g_{26}$,
\newline
\newline
$3\cdot 2^{11} \sum s(K_{1,3},G_v^-)=n^{\underline 5}-4(n-2)^{\underline 3}g_{2}+6(n-3)^{\underline 2}g_{3}-4(n-3)^{\underline 2}g_{4}-3(n-4)g_{6}+12(n-4)g_{8}+3(n-4)g_{9}-12(n-4)g_{10}+3(n-4)g_{11}-3g_{12}+12g_{13}-12g_{14}+12g_{15}-12g_{16}-12g_{19}-2g_{20}+12g_{22}+2g_{24}+4g_{25}-3g_{26}+12g_{28}-12g_{30}-6g_{31}+3g_{32}+4g_{33}-g_{34}$,
\newline
\newline
$2^{13} \sum s(C_4,G_v^-)=n^{\underline 5}-2(n-2)^{\underline 3}g_{2}-2(n-3)^{\underline 2}g_{3}-2(n-3)^{\underline 2}g_{4}-4(n-4)g_{5}-(n-4)g_{6}+20(n-4)g_{7}-4(n-4)g_{8}-(n-4)g_{9}-6(n-4)g_{10}+5(n-4)g_{11}+5g_{12}-6g_{13}-4g_{14}-6g_{15}-4g_{16}+20g_{17}-12g_{18}-4g_{19}+10g_{20}-12g_{21}-4g_{22}-4g_{23}+10g_{24}-2g_{25}-g_{26}+4g_{27}-4g_{28}+20g_{29}-6g_{30}-2g_{31}-g_{32}-2g_{33}+g_{34}$,
\newline
\newline
$2^{11} \sum s(T_{3,1},G_v^- )=n^{\underline 5}-2(n-2)^{\underline 3}g_{2}+2(n-3)^{\underline 2}g_{3}+2(n-3)^{\underline 2}g_{4}-4(n-4)g_{5}-5(n-4)g_{6}+4(n-4)g_{7}-4(n-4)g_{8}-(n-4)g_{9}+2(n-4)g_{10}+(n-4)g_{11}+g_{12}+6g_{13}+4g_{14}+2g_{15}+4g_{16}-4g_{17}-4g_{18}-4g_{19}-2g_{20}-4g_{21}+4g_{22}-4g_{23}-2g_{24}+2g_{25}-g_{26}-4g_{27}+4g_{28}+4g_{29}+6g_{30}+2g_{31}-5g_{32}-2g_{33}+g_{34}$,
\newline
\newline
$2^{12} \sum s(T_{3,2},G_v^-)=n^{\underline 5}+2(n-3)^{\underline 2}g_{3}+4(n-3)^{\underline 2}g_{4}-8(n-4)g_{5}-7(n-4)g_{6}+8(n-4)g_{7}-4(n-4)g_{8}+3(n-4)g_{9}-8(n-4)g_{10}-(n-4)g_{11}+g_{12}-4g_{13}-4g_{14}+8g_{15}+12g_{16}+4g_{19}+2g_{20}+4g_{22}+8g_{23}-2g_{24}-4g_{25}-3g_{26}-12g_{28}-8g_{29}+4g_{30}-2g_{31}+7g_{32}-g_{34}$,
\newline
\newline
$2^{7} \sum s(P_3,G_v^-)=n^{\underline 4}-2(n-2)^{\underline 2}g_{2}-g_{6}+4g_{7}-g_{9}-2g_{10}+g_{11}$,
\newline
\newline
$2^{8} \sum s(K_1,G_v^+)s(P_3,G_v^-)=n^{\underline 5}-(n-2)^{\underline 3}g_{2}-3(n-3)^{\underline 2}g_{3}+3(n-4)g_{5}+2(n-4)g_{7}+(n-4)g_{8}-(n-4)g_{9}-2(n-4)g_{10}+(n-4)g_{11}-g_{12}-2g_{13}+3g_{14}-g_{15}+2g_{16}-2g_{17}-g_{18}+g_{19}-g_{20}+3g_{21}-3g_{22}+g_{23}$,
\newline
\newline
$2^{5} \sum s(K_1,G_v^+)s(K_2,G_v^-)=n^{\underline 4}-3(n-3)g_{3}+(n-3)g_{4}+g_{5}+g_{6}-2g_{7}+g_{8}$,
\newline
\newline
$2^{6} \sum s(K_1,G_v^+)^2s(K_2,G_v^-)=n^{\underline 4}(n-2)+(n-2)^{\underline 3}g_{2}-2(n-3)(2n-5)g_{3}+(n-2)^{\underline 2}g_{4}+2g_{5}+2(n-3)g_{6}-4(n-3)g_{7}+2(n-3)g_{8}+g_{12}+g_{13}-2g_{14}+g_{15}$.

\end{proof}

\end{document}